# Reducibility mod $p$ of integral closed subschemes in projective spaces – an application of arithmetic Bézout

Reinie Erné

March 23, 2000

**Abstract.** In [4], we showed that we can improve results by Emmy Noether and Alexander Ostrowski ([8]) concerning the reducibility modulo $p$ of absolutely irreducible polynomials with integer coefficients by giving the problem a geometric turn and using an arithmetic Bézout theorem ([2]). This paper is a generalization of [4], where we show that combining the methods of [4] with the theory of Chow forms leads to similar results for flat, equidimensional, integral, closed subschemes of arbitrary codimension in $\mathbf{P}^s_{\mathbf{Z}}$.

**Introduction.** Let $K$ be a number field with ring of integers $R$, and $Z$ a flat, equidimensional, integral, closed subscheme of dimension $r+1$ and degree $d$ in $\mathbf{P}^s_R$ ($s,d \geq 2$), with absolutely irreducible generic fiber. One can show that the fiber $Z_{k(p)}$ is also absolutely irreducible for all but finitely many prime ideals $p$ of $R$ (e.g. [5, Theorem 9.7.7] and [6, Theorem 4.10]).

We would like to bound the (product of the) norms of the prime ideals $p$ of $R$ for which the fiber $Z_{k(p)}$ is not absolutely irreducible in terms of the projective height of $Z$, as defined in [2]. In this paper, using arithmetic intersection theory, we solve for any fixed $n < d$ the analogous problem obtained by replacing "absolutely irreducible" by "is not a union of two closed subschemes of degrees $n$ and $d - n$, respectively". To prove this theorem, we use Chow forms, and translate the problem to bounding the height of an intersection in some projective space. Thus, the proof becomes a straightforward application of an arithmetic Bézout Theorem for non-proper intersections given in [2], which reduces it to bounding degrees and heights of specific cycles in terms of the data provided.

**Thanks.** I would like to thank Hervé Billard for asking me about higher codimensions when I gave a talk in Brest on the irreducibility of hypersurfaces. I would also like to thank Teresa De Diego, who took the time during a conference in Obernai to tell me about Chow forms.

**Some notation.** Given a ring $R$ as above, and a locally free $\mathcal{O}_{\operatorname{Spec}(R)}$-module $\mathcal{E}$ of finite rank $s+1$ ($s \geq 0$), let $\mathbf{P}(\mathcal{E}) = \mathbf{Proj}_{\operatorname{Spec}(R)}(\operatorname{Sym}(\mathcal{E}^\vee))$ be the associated *space of lines*, where $\mathcal{E}^\vee$ denotes the dual sheaf of $\mathcal{E}$, and let $\pi$ denote its structural morphism. We suppose $\mathcal{E}$ endowed with a Hermitian metric $h$, and endow $\mathcal{E}^\vee$ with the dual metric.

Let $r$ be a positive integer. For $i = 0, \ldots, r$, let $\mathbf{P}_i = \mathbf{P}(\mathcal{E})$, and $\mathbf{P}_i^\vee = \mathbf{P}(\mathcal{E}^\vee)$. We endow the canonical quotient line bundle $\mathcal{O}(1)$ on $\mathbf{P}_i^\vee$ with the quotient metric (cf [2, 3.1.2.3]), and let $\overline{M_i}$ be the pullback of the resulting Hermitian line bundle $\overline{\mathcal{O}(1)}$ on $\mathbf{P}_i^\vee$ to $\prod_{i=o}^{r} \mathbf{P}_i^\vee$.

Finally, for $x \in \mathbf{N}$, let $F_{x,r}(\mathcal{E}) := \otimes_{i=0}^r \operatorname{Sym}^x(\mathcal{E})$

**Chow divisors and forms.** By [2, 4.3.1], we can associate to each non-zero algebraic cycle $Z \in Z_{r+1}(\mathbf{P}(\mathcal{E}))$ a *Chow divisor* $\operatorname{Ch}_{\mathbf{1}}(Z)$ (where $\mathbf{1} = (1, \ldots, 1) \in \mathbf{Z}^{r+1}$) in $Z^1(\prod_{i=0}^r \mathbf{P}_i^\vee)$, which is effective (resp. flat, resp. flat and irreducible) if such is the case for $Z$.

Let $Z$ now be a non-zero effective cycle of degree $x$ in $Z_{r+1}(\mathbf{P}(\mathcal{E}))$. Generically, the associated Chow divisor $\operatorname{Ch}_{\mathbf{1}}(Z)_K$ is the divisor of a non-zero multihomogeneous form $\phi_{\mathbf{1},Z_K}$ in

$$H^0\left(\prod_{i=0}^r (\mathbf{P}_i^\vee)_K, \otimes M_{i,K}^x\right) \cong F_{x,r}(\mathcal{E})_K,$$

called the *Chow form* of $Z_K$. Thus we can associate a point of $\mathbf{P}(F_{x,r}(\mathcal{E}))(K)$ to each non-zero effective cycle of degree $x$ in $Z_{r+1}(\mathbf{P}(\mathcal{E}))$. If the class number of $K$ is one, there exists a generalized Chow form $\phi_{\mathbf{1},Z}$ over $R$, for which $\operatorname{Ch}_{\mathbf{1}}(Z) = \operatorname{div}(\phi_{\mathbf{1},Z})$ in $\prod \mathbf{P}_i^\vee$. Similarly, for every point $t$ of $\operatorname{Spec}(R)$, we can define Chow divisors and forms for the cycles contained in the fibre above $t$ ([2, 4.3.2]).

If $Z$ is moreover flat over $\operatorname{Spec}(R)$, we have the following result:

**Proposition.** *Let $Z \in Z_{r+1}(\mathbf{P}(\mathcal{E}))$ be a flat, integral, closed subscheme of $\mathbf{P}(\mathcal{E})$ of degree $x$, with Chow divisor $\operatorname{Ch}_{\mathbf{1}}(Z)$. Let $\phi_K$ be the Chow form of $Z_K$. Let $[\phi_K] \in \mathbf{P}(F_{x,r}(\mathcal{E}))(K)$ be the corresponding point, and $P_Z$ its Zariski closure in $\mathbf{P}(F_{x,r}(\mathcal{E}))$. Then for every point $t$ of $\operatorname{Spec}(R)$, the fiber $P_{Z,t}$ is the point of $\mathbf{P}(F_{x,r}(\mathcal{E}))_t$ corresponding to the Chow form $\phi_t$ of $Z_t$.*

**Proof.** It suffices to note that by construction, we have $\operatorname{Ch}_{\mathbf{1}}(Z_t) = \operatorname{Ch}_{\mathbf{1}}(Z)_t$ for every point $t$ of $\operatorname{Spec}(R)$ ([2, 4.3.2]). In particular, as $Z$ is flat, the Zariski closure of $\operatorname{div}(\phi_K) = \operatorname{Ch}_{\mathbf{1}}(Z_K)$ is $\operatorname{Ch}_{\mathbf{1}}(Z)$. ◇

**Components of degree n.** Let $d \in \mathbf{N}_{>0}$, and fix integers $1 \leq n \leq d-1$ and $0 \leq r \leq s$. Let us simplify the notation by setting $F_x := F_{x,r}(\mathcal{E})$ for every $x$. Consider the morphism

$$\psi : \mathbf{P}(F_n) \times \mathbf{P}(F_{d-n}) \to \mathbf{P}(F_d)$$

defined by taking the product on sections (seen as multihomogeneous forms on $\prod \mathbf{P}_i^\vee$). Let $\mathcal{W}$ denote the image of $\psi$.

Let $Z$ be a flat, integral, closed subscheme of degree $d$ in $Z_{r+1}(\mathbf{P}(\mathcal{E}))$, and let $P_Z$ be as in the proposition. By dimension arguments and the proposition, the intersection of $P_Z$ and $\mathcal{W}$ is either $P_Z$, if $Z_{\overline{K}}$ has a component (irreducible or not) of degree $n$, or a finite number of closed points whose images under the structural morphism $\pi : \mathbf{P}(F_d) \to \operatorname{Spec}(R)$ are the prime ideals $q_1, \ldots, q_v$ above which the fiber of $Z \to \operatorname{Spec}(R)$ has such a component.

Before stating the theorem, we note that if $\mathcal{E}$ is isomorphic to $R^{s+1}$, then each vector bundle $F_x$ is free, and can be endowed, in a natural way, with a basis $\mathcal{B}_x$ ([2, p. 985]). Indeed, in this case, $F_x$ is a space of multihomogeneous forms as described in [2, 4.3.13], whose basis is formed by the monomials. We will use this basis to identify $F_x$ with $R^{N_x+1}$ (where $N_x := \operatorname{rk}(F_x) - 1$).

The following theorem only deals with the trivial vector bundle, i.e. $\mathcal{E} = R^{s+1}$, endowed with the standard Hermitian metric.

**Theorem.** *Let $Z \in Z_{r+1}(\mathbf{P}_R^s)$ be a flat, integral, closed subscheme of $\mathbf{P}_R^s$ of dimension $r+1$ and degree $d$ ($s, d \geq 2, r \geq 0$), and $n \in \{1, \ldots, d-1\}$ an integer such that $Z_{\overline{K}}$ cannot be written as the union of two closed subschemes of degrees $n$ and $d-n$, respectively. Let $q_1, \ldots, q_v$ be the distinct prime ideals of $R$ above which the geometric fiber of $Z$ can be written as such a union. Setting $N_{x,r,s} := \operatorname{rk}(\otimes_{i=0}^r \operatorname{Sym}^x(R^{s+1})) - 1$, we have*

$$(1) \quad \log \prod_{j=1}^v N(q_j) \leq \frac{1}{1+\delta_{n,d-n}} \binom{N_{n,r,s} + N_{d-n,r,s}}{N_{n,r,s}} h_K(Z) + O(1)$$

*when $h_K(Z)$ tends to infinity, where $h_K$ is the projective height associated to the standard Hermitian metric on $R^{s+1}$, as defined in [2, 4.1.1] (see also [3, 2.1.5]), and $\delta$ is the Kronecker delta function. Moreover, we can replace the $O(1)$ by an explicit function of $s, d, r$, and $n$ (see the proof).*

**Remark.** For the hypersurface case ($r = s$), we find a stricter bound in [4], due to the fact that horizontal hypersurfaces (which correspond to the flat integral closed subschemes here) are (directly) parametrized by a projective space, making it unnecessary to use Chow forms. The $M_x$ used there correspond to the $N_{x,r,s}$ for $r = 0$ in this paper.

**Proof.** As noted before, the set $\{q_1, \ldots, q_v\}$ is the support of $\pi(P_Z \cap \mathcal{W})$ in $\mathrm{Spec}\,(R)$. In particular, $\log \prod N(q_i) = h_K(|P_Z \cap \mathcal{W}|)$. By the arithmetic Bézout theorem [2, 5.5.1.iii], we have

$$h_K(|P_Z \cap \mathcal{W}|) \leq \deg_K(P_Z) h_K(\mathcal{W}) + h_K(P_Z) \deg_K(\mathcal{W})$$
$$+ \frac{1}{2}[K:\mathbf{Q}] \deg_K(P_Z) \deg_K(\mathcal{W})(M_d + 1)\log(2).$$

By definition of $P_Z$, its degree equals one. Using the further shortened, and somewhat misleading, notation $N_x := N_{x,r,s}$, we find that the other terms on the right can be bounded as follows:

(1) $h_K(P_Z) = h_K(Z) + d\,[K:\mathbf{Q}](\sigma_r + (r+1)\log(s+1))$,
where $\sigma_x = (1/2)(x+1)\sum_{m=2}^{x+1}(1/m)$.

(2) $\deg_K(\mathcal{W}) = \frac{1}{1+\delta_{n,d-n}}\binom{N_n+N_{d-n}}{N_n}$

(3) $h_K(\mathcal{W}) \leq \frac{[K:\mathbf{Q}]}{1+\delta_{n,d-n}} \frac{N_n+N_{d-n}+1}{2} \binom{N_n+N_{d-n}}{N_n} \log\left((d+1)^{3(r+1)(s+1)} \frac{(N_n+1)(N_{d-n}+1)}{N_n+N_{d-n}+1}\right)$,

leading to the result of the theorem. $\diamond$

**Proof of (1).** Let $\{a_I\}$ be the coefficients of $P_{Z,K}$ (i.e. of the form $\phi_{\mathbf{1},Z_K}$) in the basis $\mathcal{B}_d$. Then

$$h_K(P_Z) = \sum_\sigma \log\left(\sum |a_I|^2\right)^{1/2} - \sum_p \min_I v_p(a_I) \log N(p).$$

Another height associated to $\mathcal{B}_d$ ([2, 4.3.4.1]) is

$$h_\mathcal{B}(P_Z) := h_\mathcal{B}(\mathrm{Ch}_{\mathbf{1}}(Z)) = \sum_\sigma \log\left(\sum |a_I|\right) - \sum_p \min_I v_p(a_I) \log N(p).$$

Clearly, we have $h_K(P_Z) \leq h_\mathcal{B}(P_Z)$. By [2, Theorem 4.3.8, (4.3.33), and (4.1.2)],

$$h_\mathcal{B}(P_Z) \leq h_K(Z) + d\,[K:\mathbf{Q}](\sigma_r + (r+1)\log(s+1)).$$

In particular, we have $h_K(P_Z) = h_K(Z) + \mathcal{O}(1)$.

**Remark.** Before giving the proofs of (2) and (3), let us note that the morphism $\psi$ was used under the notation $\phi_n$ in [4], where the degree and height of its image were bounded explicitely. Here we give only sketches of the proofs of (2) and (3), the details can be found in [loc.cit.].

**Proof of (2).** Let $f_n$ (resp. $f_{d-n}$) denote the projection from $\mathbf{P}(F_n) \times \mathbf{P}(F_{d-n})$ onto the first (resp. second) coordinate. Using intersection theory, we find

$$\deg(\psi) \deg_K(\mathcal{W}) = \deg \left( c_1 \mathcal{O}_{\mathbf{P}(F_d)}(1)^{N_n + N_{d-n}} \cdot [\psi_*(\mathbf{P}(F_n) \times \mathbf{P}(F_{d-n}))] \right),$$

where $c_1 \mathcal{O}_{\mathbf{P}(F_d)}(1)$ is the first Chern class of $\mathcal{O}_{\mathbf{P}(F_d)}(1)$. By the projection formula, and the fact that $\psi^* \mathcal{O}_{\mathbf{P}(F_d)}(1) = f_n^* \mathcal{O}_{\mathbf{P}(F_n)}(1) \otimes f_{d-n}^* \mathcal{O}_{\mathbf{P}(F_{d-n})}(1)$, this implies that

$$\deg_K(\mathcal{W}) = \frac{1}{1 + \delta_{n,d-n}} \binom{N_n + N_{d-n}}{N_n}.$$

**Proof of (3).** As in the proof of (2), we use intersection theory, but this time with metrics. By [2, 4.1.2 and Proposition 2.3.1], we have

$$h_K(\mathcal{W}) = [K : \mathbf{Q}] \sigma_{N_n + N_{d-n}} \deg_K(\mathcal{W})$$
$$+ \frac{1}{\deg(\psi)} \widehat{\deg} \left( \widehat{c}_1 \left( \psi^* \overline{\mathcal{O}_{\mathbf{P}(F_d)}(1)} \right)^{N_n + N_{d-n} + 1} \mid \mathbf{P}(F_n) \times \mathbf{P}(F_{d-n}) \right).$$

The arithmetic degree on the right is the projective height of $\mathbf{P}(F_n) \times \mathbf{P}(F_{d-n})$ associated to the line bundle $\mathcal{L} := \psi^* \mathcal{O}_{\mathbf{P}(F_d)}(1)$ endowed with the pullback $\rho$ under $\psi$ of the standard Hermitian metric on $\mathcal{O}_{\mathbf{P}(F_d)}(1)$. We will bound this height in two steps, using a comparison of metrics on $\mathcal{L}$. First, let $(\mathcal{L}, \rho')$ denote the line bundle $\mathcal{L}$ endowed with the product metric obtained by taking the standard Hermitian metrics on $\mathcal{O}_{\mathbf{P}(F_n)}(1)$ and $\mathcal{O}_{\mathbf{P}(F_{d-n})}(1)$. The associated projective height is

$$h_{(\mathcal{L}, \rho')}(\mathbf{P}(F_n) \times \mathbf{P}(F_{d-n})) := \widehat{\deg} \left( \widehat{c}_1(\mathcal{L}, \rho')^{N_n + N_{d-n} + 1} \mid \mathbf{P}(F_n) \times \mathbf{P}(F_{d-n}) \right),$$

which, by the projection formula and the decomposition of the metrized line bundle $(\mathcal{L}, \rho')$ as a (tensor)product, equals

$$[K : \mathbf{Q}] \left( \binom{N_n + N_{d-n}}{N_n} \sigma_{N_{d-n}} + \binom{N_n + N_{d-n}}{N_n} \sigma_{N_n} \right).$$

The second step in bounding the height that we want consists of comparing the norms $\|\cdot\|$ and $\|\cdot\|'$ associated to $\rho$, resp. $\rho'$. Let $\varphi : (\mathbf{P}(F_n) \times \mathbf{P}(F_{d-n}))(\mathbf{C}) \to \mathbf{R}$ be defined by $(\|\cdot\|')^2 = \exp(\varphi) \|\cdot\|^2$. For each embedding $\sigma : K \hookrightarrow \mathbf{C}$, and $(a, b) = ((a_0 : \ldots : a_{N_n}), (b_0 : \ldots : b_{N_{d-n}}))$ in $(\mathbf{P}(F_n) \times \mathbf{P}(F_{d-n}))_\sigma(\mathbf{C})$, we let $f_a$, resp. $g_b$, be the corresponding multihomogeneous polynomial (in $(r+1)(s+1)$ variables). We have

$$\exp(\varphi_\sigma(a,b)) = \left(\frac{L_2(f_a g_b)}{L_2(f_a) L_2(g_b)}\right)^2,$$

where the $L_2$-norm $L_2(f)$ of a (multi)homogeneous polynomial $f = \sum c_I X^I$ is $(\sum c_I \overline{c_I})^{1/2}$. From results of [7, 3.2], we can now deduce that

$$\sup_{(a,b)}(\varphi_\sigma(a,b)) \leq 3(r+1)(s+1)\log(d+1)$$

for every $\sigma : K \hookrightarrow \mathbf{C}$. The last step consists of combining this inequality with the results of [2, Proposition 3.2.2] and [1, Lemma 2.6.ii] (see also [4]) to obtain

$$h_{(\mathcal{L},\rho)}(\mathbf{P}(F_n) \times \mathbf{P}(F_{d-n})) \leq h_{(\mathcal{L},\rho')}(\mathbf{P}(F_n) \times \mathbf{P}(F_{d-n}))$$
$$+ [K:\mathbf{Q}]\frac{N_n + N_{d-n} + 1}{2}\deg(\psi)\deg_K(\mathcal{W})\,3\,(r+1)(s+1)\log(d+1),$$

which, after some simplification, leads to the bound for $h_K(\mathcal{W})$ stated in (3).

Reinie Erné
Institut Mathémathiques
Université de Rennes I
Campus de Beaulieu
35042 Rennes cedex
France
email: erne@univ-rennes1.fr